\documentclass[a4paper,11pt]{article}
\usepackage[top=2.5cm,bottom=2.5cm,left=2.5cm,right=2.5cm]{geometry}
\usepackage{amsfonts}
\usepackage{latexsym,amsmath,amssymb,amsfonts,epsfig,psfrag,url,graphics,ifpdf,multicol}
\usepackage{color,colordvi,amscd,amsthm} 
\usepackage{tikz}
\usepackage[numbers,sort&compress]{natbib}
\usepackage{indentfirst,graphics,epsfig}
\usepackage{graphicx}
\usepackage{graphics}
\usepackage{graphicx,psfrag}
\usepackage{graphics}

\usepackage{tabularx}
\usepackage{booktabs}
\usepackage{floatrow}
\usepackage{multirow}
\usepackage{graphicx}
\usepackage{caption}
\usepackage[ruled]{algorithm2e}
\newfloatcommand{capbtabbox}{table}[][\FBwidth]

\newtheorem{theo}{Theorem}[section]
\newtheorem{lem}[theo]{Lemma}
\newtheorem{coro}[theo]{Corollary}

 \theoremstyle{definition}
\newtheorem{defn}[theo]{Definition}

\theoremstyle{remark}

\def\pf{\noindent {\bf Proof.} }

\newcounter{casenum}[theo]

\newcounter{subcasenum}[theo]

\newcounter{claimnum}[theo]

\allowdisplaybreaks
\usepackage{indentfirst,subfig}

\begin{document}
	\thispagestyle{empty}
	\captionsetup[figure]{labelfont={bf},name={Fig.},labelsep=period}

\begin{center} {\large\sc
The maximum number of maximum dissociation sets in trees}
\end{center}
\begin{center}
{
  {\small Jianhua Tu$^a$, Zhipeng Zhang$^a$, Yongtang Shi$^{b}$\footnote{Corresponding author.\\\indent \ \  E-mail: tujh81@163.com (J. Tu); 18622423315@163.com (Z. Zhang); shi@nankai.edu.cn (Y. Shi)}}\\[2mm]

{\small $^a$ Department of mathematics, Beijing University of Chemical Technology \\
\hspace*{1pt} Beijing 100029, China\\
$^b$ Center for Combinatorics and LPMC\\ Nankai University, Tianjin, P.R. China 300071
} \\[2mm]
}
\end{center}

\begin{center}

\begin{abstract}
A subset of vertices is a {\it maximum independent set} if no two of the vertices are adjacent and the subset has maximum cardinality. A subset of vertices is called a {\it maximum dissociation set} if it induces a subgraph with vertex degree at most 1, and the subset has maximum cardinality. Zito [J. Graph Theory {\bf 15} (1991) 207--221] proved that the maximum number of maximum independent sets of a tree of order $n$ is $2^{\frac{n-3}{2}}$ if $n$ is odd, and $2^{\frac{n-2}{2}}+1$ if $n$ is even and also characterized all extremal trees with the most maximum independent sets, which solved a question posed by Wilf.
Inspired by the results of Zito, in this paper, by establishing four structure theorems and a result of $k$-K\"{o}nig-Egerv\'{a}ry graph, we show that the maximum number of maximum dissociation sets in a tree of order $n$ is

\begin{center}
$\left\{
  \begin{array}{ll}
    3^{\frac{n}{3}-1}+\frac{n}{3}+1, & \hbox{if $n\equiv0\pmod{3}$;} \\
    3^{\frac{n-1}{3}-1}+1, & \hbox{if $n\equiv1\pmod{3}$;} \\
    3^{\frac{n-2}{3}-1}, & \hbox{if $n\equiv2\pmod{3}$,}
  \end{array}
\right.$
\end{center}
and also give complete structural descriptions of all extremal trees on which these maxima are achieved.

\vspace{5mm}

\end{abstract}
\end{center}

\baselineskip=0.24in

\section{Introduction}

In this paper, we consider undirected labeled graphs without loops or multiple edges and use standard graph-theoretic terminology (see \cite{Bondy2008}). An {\it independent set} of a graph $G$ is a set of vertices no two of which are joined by an edge. An independent set is called {\it maximal} if it cannot be contained in any other independent set, and is {\it maximum} if it has maximum cardinality.  The {\it independence number $\alpha(G)$} of $G$ is the cardinality of a maximum independent set of $G$.

A subset of vertices in a graph $G$ is called a {\it dissociation set} if it induces a subgraph with vertex degree at most 1.
The {\it dissociation number} of a graph $G$, denoted by $diss(G)$, is the cardinality of a maximum dissociation set of $G$.
The problem of computing $diss(G)$ (dissociation number problem) has been introduced by Yannakakis \cite{Yannakakis1981} and was shown to be NP-complete for the class of bipartite graphs.  Actually, Orlovich, Dolgui, Finke, Gordon and Werner \cite{Orlovich2011} showed that it remains NP-hard even in planar line graphs of planar bipartite graphs. Cameron and Hell \cite{Cameron2006} proved that the problem can be solved in polynomial time for chordal graphs, weakly
chordal graphs, asteroidal triple-free graphs and interval-filament graphs. The complexity of the problem on some classes of graphs has been studied in \cite{Alekseev2007,Boliac2004,Bresar2011,Cameron2006,Orlovich2011,Papadimitriou}. Note that a set $S$ of vertices of a graph $G$ is a dissociation set if and only if its complement $V(G)\setminus S$ is a so-called {\it 3-path vertex cover}, that is, a set of vertices intersecting every path of order 3 in $G$. The 3-path vertex cover problem is to find a minimum 3-path vertex cover in a given graph and has been well studied \cite{BaiTuShi,Bresar2011,Kardos2011,Xiao2017}.

In 1986, Wilf \cite{Wilf1986} determined the maximum number of maximal independent sets in a tree. Later, Sagan \cite{Sagan1988} gave a short proof and characterized the extremal trees. F\"uredi \cite{Furedi1987} found the maximum number of maximal independent sets for connected graphs on $n>50$ vertices. Independently, Griggs, Grinstead and Guichard \cite{Griggs1988} determined the maximum number of maximal independent sets for connected graphs on $n$ vertices for all values of $n$ and completely characterized all the extremal graphs. For more results on the maximum number of maximal independent sets, we refer to \cite{Koh2008,Law,Liu,Sagan2006,Wagner,Wloch2008}.
Zito \cite{Zito1991} proved that the maximum number of maximum independent sets of a tree of order $n$ is $2^{\frac{n-3}{2}}$ if $n$ is odd, and $2^{\frac{n-2}{2}}+1$ if $n$ is even and she also characterized all extremal trees with the most maximum independent sets, which solved a question posed by Wilf \cite{Wilf1986}. Alvarado, Dantas, Mohr and Rautenbach \cite{Alvarado2019} showed that every tree with independence number $\alpha$ has at most $2^{\alpha-1}+1$ maximum independent sets.

Inspired by the results of Zito \cite{Zito1991},
in this paper, we consider the analogous problem of finding the maximum number of maximum dissociation sets and the extremal graphs for trees of order $n$. By establishing four structure theorems and a result of $k$-K\"{o}nig-Egerv\'{a}ry graph, we show that the maximum number of maximum dissociation sets in a tree of order $n$ is
\begin{center}
$\left\{
  \begin{array}{ll}
    3^{\frac{n}{3}-1}+\frac{n}{3}+1, & \hbox{if $n\equiv0\pmod{3}$;} \\
    3^{\frac{n-1}{3}-1}+1, & \hbox{if $n\equiv1\pmod{3}$;} \\
    3^{\frac{n-2}{3}-1}, & \hbox{if $n\equiv2\pmod{3}$,}
  \end{array}
\right.$
\end{center}
and also characterize the structure of the extremal trees with the most maximum dissociation sets.

The paper is organized as follows. In next section, we introduce and study the $k$-K\"{o}nig-Egerv\'{a}ry graphs. We show that any forest is a $k$-K\"{o}nig-Egerv\'{a}ry graph, which plays a key role in presenting structure theorems in Section 3. In Section 3, four structure theorems are established.
In Section 4, we apply these structure theorems to find the families of trees with the most maximum dissociation sets.

\section{$k$-K\"{o}nig-Egerv\'{a}ry graphs}

Let $G$ be a graph. The set of neighbors of a vertex $v$ in $G$ is denoted by $N_G(v)$. Let $U$ be a set of vertices in $G$. The set of all neighbors of the vertices in $U$ is denoted by $N_G(U)$. For a positive integer $k$, a {\it $k$-path} is a (not necessarily induced) path of order $k$. Let $P_k$ be the path with $k$ vertices.

A {\it matching} in a graph $G$ is a set of edges no two of which share one common vertex. The {\it matching number $\mu(G)$} is the cardinality of a maximum matching of $G$. The famous K\"{o}nig-Egerv\'{a}ry theorem states that for any bipartite graph $G$, the sum of independence number $\alpha(G)$ and matching number $\mu(G)$ equals $|V(G)|$. A graph $G$ is called a {\it K\"{o}nig-Egerv\'{a}ry graph} if $\alpha(G)+\mu(G)=|V(G)|$.
Clearly, every bipartite graph is a K\"{o}nig-Egerv\'{a}ry graph.  K\"{o}nig-Egerv\'{a}ry graphs have been extensively studied \cite{Cardoso2017,Jarden2017,Levit2011,Lovasz1983,Sterboul1979}.

A {\it $k$-matching} in a graph $G$ is a set of vertex-disjoint $k$-paths in $G$, and the {\it $k$-matching number $\mu_k(G)$} of $G$ is the cardinality of a maximum $k$-matching in $G$. A {\it $k$-vertex cover} is a set of vertices of $G$ intersecting every $k$-path, and the {\it $k$-vertex cover number} $\tau_k(G)$ of $G$ is the cardinality of a minimum $k$-vertex cover in $G$. A {\it $k$-independent set} is a set $S$ of vertices such that the subgraph induced by $S$ contains no $k$-paths, and the {\it $k$-independence number $\alpha_k(G)$} is the cardinality of a maximum $k$-independent set in $G$. Note that $\mu_2(G)$, $\alpha_2(G)$ and $\alpha_3(G)$ are exactly the matching number $\mu(G)$, independence number $\alpha(G)$ and dissociation number $diss(G)$, respectively.

It can be easily seen that for any graph $G$, $\alpha_k(G)+\tau_k(G)=|V(G)|$, $\mu_k(G)\leq \tau_k(G)$ and $\alpha_k(G)+\mu_k(G)\leq|V(G)|$.
Now we introduce a generalization of the K\"{o}nig-Egerv\'{a}ry graphs, which are called $k$-K\"{o}nig-Egerv\'{a}ry graphs.

\begin{defn}\label{def1}
Given a positive integer $k\geq 2$, a graph $G$ is called a {\it $k$-K\"{o}nig-Egerv\'{a}ry graph} if $\alpha_k(G)+\mu_k(G)=|V(G)|$.
\end{defn}

Clearly, a graph $G$ is $k$-K\"{o}nig-Egerv\'{a}ry if and only if $\tau_k(G)=\mu_k(G)$. In a rooted tree, the {\it level} of a vertex $v$ is the length of the unique path from the root to the vertex $v$, denoted by $\ell(v)$.

\begin{theo}\label{2.1}
For a positive integer $k\geq 2$, any forest is a $k$-K\"{o}nig-Egerv\'{a}ry graph.
\end{theo}

\pf It suffices to show that any tree $T$ is a $k$-K\"{o}nig-Egerv\'{a}ry graph. We prove this by induction on the number of vertices of $T$.

If $|V(T)|\leq k-1$, then $\tau_k(T)=\mu_k(T)=0$ and $T$ is a $k$-K\"{o}nig-Egerv\'{a}ry graph.

Assume that the result is true for all trees with fewer than $n$ vertices. Let $T$ be a tree with $n$ ($n\geq k$) vertices. Change the tree $T$ into a rooted tree by choosing any vertex as the root. Suppose that a vertex $u$ is chosen such that there is a $k$-path in the subtree $T_u$ rooted at $u$ and, subject to this condition, the level of $u$ is as large as possible.
It is easy to see that any $k$-path in the subtree $T_u$ must contain the vertex $u$. Let $T':=T-T_u$.

Suppose that $F$ is a minimum $k$-vertex cover of $T$. Then $F\cap V(T_u)\neq \emptyset$ and $F_1:=(F\setminus V(T_u))\cup\{u\}$ is also a minimum $k$-vertex cover of $T$. Thus $F_1\setminus\{u\}$ is a $k$-vertex cover of $T'$ and $\tau_k(T')\leq \tau_k(T)-1$. On the other hand, if $F'$ is a minimum $k$-vertex cover of $T'$, then $F'\cup\{u\}$ is a $k$-vertex cover of $T$. Thus, $\tau_k(T)\leq \tau_k(T')+1$. So, $\tau_k(T)=\tau_k(T')+1$.

Let $\mathcal{M}$ be a maximum $k$-matching of $T$ and $P$ a $k$-path of $T_u$. We can change $\mathcal{M}$ into a maximum $k$-matching $\mathcal{M}_1$ such that $P\in \mathcal{M}_1$ and $\mathcal{M}_1\setminus\{P\}$ is a $k$-matching of $T'$. Thus, $\mu_k(T')\geq \mu_k(T)-1$. On the other hand, if $\mathcal{M}'$ is a maximum $k$-matching of $T'$, then $\mathcal{M'}\cup\{P\}$ is a $k$-matching of $T$ and $\mu_k(T)\geq \mu_(T')+1$. So, $\mu_(T)=\mu_k(T')+1$.

By the induction hypothesis, $\tau_k(T')=\mu_k(T')$. Thus $\tau_k(T)=\mu_k(T)$ and $T$ is a $k$-K\"{o}nig-Egerv\'{a}ry graph.  \qed

\section{The structure theorems}

In this section, we give four structure theorems concerning $\alpha_3$-critical edges and three types of vertices.

An edge $e$ of a graph $G$ is called {\it $\alpha_3$-critical} if $\alpha_3(G-e)>\alpha_3(G)$. An edge $e$ of a graph $G$ is called {\it $\mu_3$-critical} if $\mu_3(G-e)<\mu_3(G)$.
A subgraph of a graph $G$ is called {\it critical} if all its edges are $\alpha_3$-critical in $G$. An $\alpha_3$-critical edge of a graph $G$ is called {\it insulated} if it is not adjacent to any other $\alpha_3$-critical edge of $G$. If $M$ is a $3$-matching, an edge of a 3-path of $M$ is said to be {\it covered} by $M$, and each vertex of a 3-path of $M$ is said to be {\it saturated} by $M$.

A vertex in a graph $G$ is said to be \emph{flexible} if it is in some but not all maximum dissociation sets of $G$,
in what follows we use $\mathcal{F}_G$ to denote the set of all flexible vertices of $G$.
A vertex is called \emph{static} if it is either in all maximum dissociation sets or in no maximum dissociation sets. If a vertex is in all maximum dissociation sets call it \emph{static-included}, and we use $\mathcal{A}_G$ to denote the set of all static-included vertices of $G$.
If a vertex is in no maximum dissociation set call it \emph{static-excluded}, and we use $\mathcal{N}_G$ to denote the set of all static-excluded vertices of $G$.

\begin{lem}\label{3.1}
Let $G$ be a graph and $uv$ an edge of $G$. If $uv$ is an $\alpha_3$-critical edge of $G$, then every maximum dissociation set of $G-uv$ contains both $u$ and $v$, and $\alpha_3(G-uv)=\alpha_3(G)+1$.
\end{lem}
\pf Let $F'$ be a maximum dissociation set of $G-uv$. Suppose, to the contrary, that either $u$ or $v$ is not in $F'$. Then $F'$ is also a dissociation set of $G$. Thus, $\alpha_{3}(G)\geq \alpha_{3}(G-uv)$, a contradiction. So $F'$ contains both $u$ and $v$. On the other hand, $F'-u$ is a dissociation set of $G$. Thus $\alpha_3(G)\geq \alpha_3(G-uv)-1$ and $\alpha_{3}(G-uv)=\alpha_{3}(G)+1$.\qed

Similarly, we have

\begin{lem}\label{3.2}
Let $G$ be a graph and $uv$ an edge of $G$. If $uv$ is a $\mu_3$-critical edge of $G$, then it is covered by every maximum $3$-matching of $G$, and $\mu_3(G-uv)=\mu_3(G)-1$.
\end{lem}


\begin{lem}\label{3.3}
Let $G$ be a $3$-K\"{o}nig-Egerv\'{a}ry graph and $e$ an edge of $G$. Then $e$ is $\alpha_3$-critical in $G$ if and only if:

(i) $e$ is $\mu_3$-critical in $G$,

(ii) $G-e$ is a $3$-K\"{o}nig-Egerv\'{a}ry graph.
\end{lem}

\pf Suppose that $e$ is $\alpha_3$-critical in $G$.
By Lemma \ref{3.1}, $\alpha_3(G)=\alpha_3(G-e)-1$. We have
\begin{eqnarray*}
  &&\alpha_3(G)+\mu_3(G)=\alpha_3(G-e)-1+\mu_3(G)=|V(G)|\ \ and\\
&&\alpha_3(G-e)+\mu_3(G-e)\leq |V(G-e)|=|V(G)|.
\end{eqnarray*}
Thus, $\mu_3(G)\geq \mu_3(G-e)+1$. So $e$ is $\mu_3$-critical in $G$ and $\mu_3(G)= \mu_3(G-e)+1$. Moreover, $\alpha_3(G-e)+\mu_3(G-e)=|V(G-e)|$ and $G-e$ is a $3$-K\"{o}nig-Egerv\'{a}ry graph.

On the other hand, suppose that $e$ is $\mu_3$-critical in $G$ and $G-e$ is a $3$-K\"{o}nig-Egerv\'{a}ry graph. Then $\alpha_3(G)+\mu_3(G)= |V(G)|=|V(G-e)|=\alpha_3(G-e)+\mu_3(G-e)=\alpha_3(G-e)+\mu_3(G)-1$. So, $\alpha_3(G-e)=\alpha_3(G)+1$ and $e$ is $\alpha_3$-critical in $G$. \qed

By Theorem \ref{2.1} and Lemma \ref{3.3}, we have the following corollary.

\begin{coro}\label{3.4}
Let $T$ be a tree. An edge of $T$ is $\alpha_3$-critical if and only if it is $\mu_3$-critical.
\end{coro}

The first structure theorem concerns the relationship between $\alpha_3$-critical edges and flexible vertices in trees.

\begin{theo}\label{3.5}
Let $T$ be a tree. Then

(1) every maximum dissociation set of $T$ contains at least one end-vertex of each ${\alpha _3}$-critical edge;

(2) a vertex of $T$ is flexible if and only if it is an end-vertex of an ${\alpha _3}$-critical edge.
\end{theo}

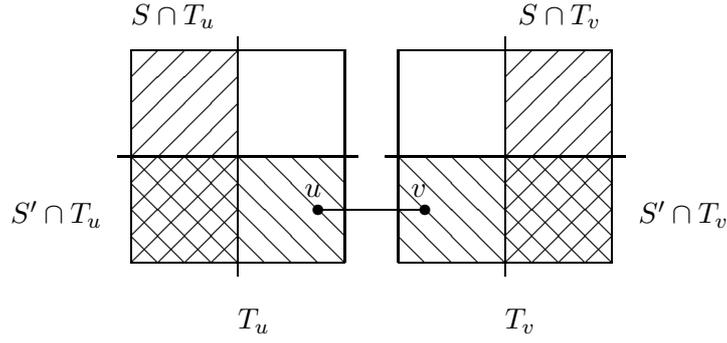
\begin{figure}[h]
	\begin{center}
		\begin{picture}(400,120)
		\multiput(180,40)(40,0){2}{\circle*{3.5}}
		\put(190,20){\line(0,1){80}}
		\put(210,20){\line(0,1){80}}
		\put(210,20){\line(1,0){80}}
		\put(190,20){\line(-1,0){80}}
		\put(210,100){\line(1,0){80}}
		\put(190,100){\line(-1,0){80}}
		\put(290,100){\line(0,-1){80}}
		\put(110,100){\line(0,-1){80}}
		\put(105,60){\line(1,0){90}}
		\put(205,60){\line(1,0){90}}
		\put(150,105){\line(0,-1){90}}
		\put(250,105){\line(0,-1){90}}	
		\put(180,40){\line(1,0){40}}
		\put(110,90){\line(1,1){10}}
		\put(110,80){\line(1,1){20}}
		\put(110,70){\line(1,1){30}}
		\put(110,60){\line(1,1){40}}
		\put(110,50){\line(1,1){40}}
		\put(110,40){\line(1,1){40}}
		\put(110,30){\line(1,1){40}}
		\put(110,20){\line(1,1){40}}
		\put(120,20){\line(1,1){30}}
		\put(130,20){\line(1,1){20}}
		\put(140,20){\line(1,1){10}}
		\put(250,90){\line(1,1){10}}
		\put(250,80){\line(1,1){20}}
		\put(250,70){\line(1,1){30}}
		\put(250,60){\line(1,1){40}}
		\put(250,50){\line(1,1){40}}
		\put(250,40){\line(1,1){40}}
		\put(250,30){\line(1,1){40}}
		\put(250,20){\line(1,1){40}}
		\put(260,20){\line(1,1){30}}
		\put(270,20){\line(1,1){20}}
		\put(280,20){\line(1,1){10}}
		\put(180,60){\line(1,-1){10}}
		\put(170,60){\line(1,-1){20}}
		\put(160,60){\line(1,-1){30}}
		\put(150,60){\line(1,-1){40}}
		\put(140,60){\line(1,-1){40}}
		\put(130,60){\line(1,-1){40}}
		\put(120,60){\line(1,-1){40}}
		\put(110,60){\line(1,-1){40}}
		\put(110,50){\line(1,-1){30}}
		\put(110,40){\line(1,-1){20}}
		\put(110,30){\line(1,-1){10}}
		\put(280,60){\line(1,-1){10}}
		\put(270,60){\line(1,-1){20}}
		\put(260,60){\line(1,-1){30}}
		\put(250,60){\line(1,-1){40}}
		\put(240,60){\line(1,-1){40}}
		\put(230,60){\line(1,-1){40}}
		\put(220,60){\line(1,-1){40}}
		\put(210,60){\line(1,-1){40}}
		\put(210,50){\line(1,-1){30}}
		\put(210,40){\line(1,-1){20}}
		\put(210,30){\line(1,-1){10}}
		\put(150,-5){$T_u$}
		\put(250,-5){$T_v$}
		\put(175,45){$u$}
		\put(215,45){$v$}
		\put(300,35){$S'\cap {T_v}$}
		\put(65,35){$S'\cap {T_u}$}
		\put(110,110){$S \cap {T_u}$}
		\put(255,110){$S \cap {T_v}$}
		\caption{\label{1} Dissociation sets ${S_1}=(S \cap {T_u}) \cup (S' \cap {T_v})$ and ${S_2}=(S \cap {T_v}) \cup (S' \cap {T_u})$.}	
		\end{picture}
	\end{center}
\end{figure}

\pf (1) Let $uv$ be an ${\alpha _3}$-critical edge in $T$, and let $S$ be a maximum dissociation set of $T$.
Suppose, to the contrary, that neither $u$ nor $v$ is in $S$.
Denote by ${T_u}$ (resp.~${T_v}$) the connected component of $T-uv$ containing $u$ (resp.~$v$).
Let $S'$ be a maximum dissociation set of $T-uv$. By Lemma \ref{3.1}, $|S'|=\alpha_3(T-uv)=\alpha_3(T)+1$.
Both ${S_1}=(S \cap {T_u}) \cup (S' \cap {T_v})$ and ${S_2}=(S \cap {T_v}) \cup (S' \cap {T_u})$ are dissociation sets of $T$, as shown in Figure \ref{1}. Since $|S_1|+|S_2|=|S|+|S'|=2\alpha_3(T)+1$, one of $S_1$ and $S_2$ must contain $\alpha_3(T)+1$ vertices, a contradiction.

\begin{figure}[h]
	\begin{center}
		\begin{picture}(400,125)
		\multiput(120,35)(80,0){3}{\circle{40}}
		\multiput(200,115)(0,0){1}{\circle*{3.5}}
		\multiput(200,55)(0,0){1}{\circle*{3.5}}
		\multiput(266,49)(0,0){1}{\circle*{3.5}}
		\multiput(134,49)(0,0){1}{\circle*{3.5}}
         \multiput(155,35)(5,0){3}{\circle*{2}}
		\multiput(235,35)(5,0){3}{\circle*{2}}
		\put(200,115){\line(0,-1){60}}
		\put(200,115){\line(1,-1){66}}
		\put(200,115){\line(-1,-1){66}}
		\put(117,45){$u_1$}
		\put(197,122){$v$}
		\put(197,45){$u_j$}
		\put(270,45){$u_k$}
		\put(120,-5){$T_1$}
		\put(200,-5){$T_j$}
		\put(280,-5){$T_k$}
		\caption{\label{2}Tree components used in the proof of Theorem \ref{3.5}(2).}	
		\end{picture}
	\end{center}
\end{figure}
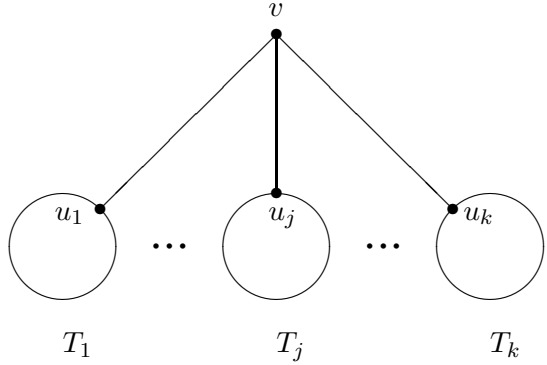

(2) Let $v$ be a flexible vertex of $T$.
Let $u_1, u_2, \ldots, u_k$ be the neighbors of $v$ in $T$, and
$T_1, T_2, \ldots, T_k$ be the connected components of $T-v$ such that $u_i \in T_i$, as shown in Figure \ref{2}.
Let $S_v$ be a maximum dissociation set of $T$ containing $v$ and $S_{\overline{v}}$ a maximum dissociation set of $T$ that does not contain $v$.
Then
\begin{eqnarray*}
 &&|\{v\}|+\sum\limits_{i=1}^k |S_v\cap T_i|=1+\sum\limits_{i=1}^k |S_v\cap T_i|= {\alpha _3}(T)\ \ and \\
&&\sum\limits_{i=1}^k |S_{\bar{v}}\cap T_i|= {\alpha _3}(T).
\end{eqnarray*}
Thus there exists $j \in \{ 1,2, \ldots ,k\} $ such that
\[\left| {{S_{\bar v}} \cap {T_j}} \right| > \left| {{S_v} \cap {T_j}} \right|.\] Thus
$S=({S_{\bar v}}\cap {T_j}) \cup ({S_v}\cap (T - {T_j}))$ is a dissociation set of $T-vu_{j}$ and $|S|>|S_{v}|=\alpha_3(T)$. It means that the edge $vu_j$ is $\alpha_3$-critical in $T$ and $v$ is incident to an $\alpha_3$-critical edge.

Let $uv$ be an ${\alpha _3}$-critical edge in $T$ and $S$ a maximum dissociation set of $T-uv$. By Lemma \ref{3.1}, both $u$ and $v$ are in $S$. It is easy to see that both $S\setminus\{u\}$ and $S\setminus\{v\}$ are maximum dissociation sets of $T$. Thus both $u$ and $v$ are flexible in $T$. \qed

The second structure theorem gives adjacency rules that determine the sets of flexible vertices, static-included vertices and $\alpha_3$-critical edges of a tree.

\begin{theo}\label{3.6}
Let $T$ be a tree. Then

(1) there is no critical $4$-path or critical $K_{1,3}$ in $T$,

(2) every end-vertex of each insulated $\alpha_3$-critical edge has exactly one neighbor in $\mathcal{A}_T$ and the neighbor must be an isolated vertex of $T[\mathcal{A}_T]$,

(3) every vertex of each critical $3$-path is not adjacent to any vertex of $\mathcal{A}_T$.

\end{theo}

\pf (1) By Corollary \ref{3.4} and Lemma \ref{3.2}, if an edge of $T$ is $\alpha_3$-critical, then it is $\mu_3$-critical and is covered by all maximum $3$-matchings of $T$. Thus, there is no critical $4$-path or critical $K_{1,3}$ in $T$.

(2) Let $uv$ be an insulated $\alpha_3$-critical edge of $T$. By Theorem \ref{3.5}(2), both $u$ and $v$ are flexible in $T$. Let $S$ be a maximum dissociation set of $T$ containing $u$. Since $S$ contains all vertices of $\mathcal{A}_T$, $u$ has at most one neighbor in $\mathcal{A}_T$ and the neighbor must be an isolated vertex of $T[\mathcal{A}_T]$. Similarly, $v$ has at most one neighbor in $\mathcal{A}_T$ and the neighbor must be an isolated vertex of $T[\mathcal{A}_T]$.

Next, we show that both $u$ and $v$ have at least one neighbor in $\mathcal{A}_T$. If the statement is not true, then we have the following two cases.

\textbf{Case 1.} Neither $u$ nor $v$ has a neighbor in $\mathcal{A}_T$.

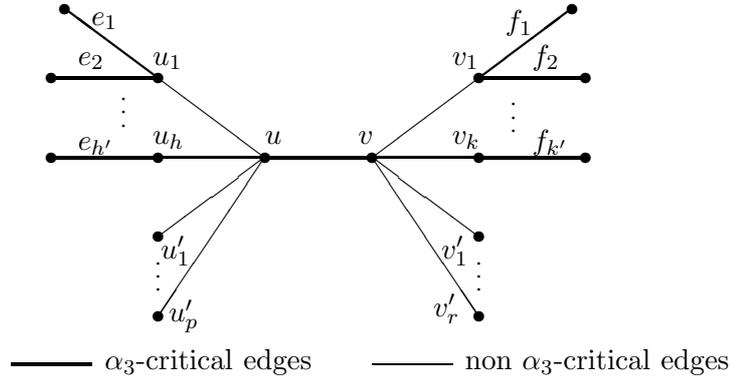
\begin{figure}[h]
	\begin{center}
		\begin{picture}(400,140)
		\thicklines
		\multiput(100,75)(40,0){6}{\circle*{3.5}}
		\multiput(100,105)(40,0){2}{\circle*{3.5}}
		\multiput(140,45)(0,-30){2}{\circle*{3.5}}
		\multiput(260,45)(0,-30){2}{\circle*{3.5}}
		\multiput(260,105)(40,0){2}{\circle*{3.5}}
		\multiput(127,87)(0,5){3}{\circle*{1}}
		\multiput(140,35)(0,-5){3}{\circle*{1}}
		\multiput(260,35)(0,-5){3}{\circle*{1}}	
		\multiput(273,85)(0,5){3}{\circle*{1}}		
		\linethickness{0.4mm}
		\put(140,105){\line(-4,3){35}}
        \put(260,105){\line(4,3){35}}
        \put(105,131){\circle*{3.5}}
        \put(295,131){\circle*{3.5}}
        \put(100,105){\line(1,0){40}}
		\put(260,105){\line(1,0){40}}
		\put(100,75){\line(1,0){40}}
		\put(180,75){\line(1,0){40}}
		\put(260,75){\line(1,0){40}}
		
		\thinlines
		\put(180,75){\line(-4,-3){40}}
		\put(180,75){\line(-2,-3){40}}
		\put(180,75){\line(-4,3){40}}
		\put(220,75){\line(4,-3){40}}
		\put(220,75){\line(2,-3){40}}
		\put(220,75){\line(4,3){40}}
		\put(140,75){\line(1,0){40}}
		\put(220,75){\line(1,0){40}}
		\put(180,80){$u$}
		\put(215,80){$v$}
        \put(270,123){$f_1$}
		\put(280,110){$f_2$}
		\put(280,78){$f_{k'}$}
        \put(115,125){$e_1$}
		\put(110,110){$e_2$}
		\put(110,78){$e_{h'}$}
		\put(250,110){$v_1$}
		\put(250,80){$v_k$}
		\put(138,80){$u_{h}$}
		\put(138,110){$u_1$}
		\put(141,37){$u'_1$}
		\put(144,13){$u'_p$}
		\put(246,37){$v'_1$}
		\put(243,15){$v'_r$}
        \linethickness{0.4mm}
        \put(85,-3){\line(1,0){30}}
		\put(120,-5){$\alpha_3$-critical edges}
        \thinlines
        \put(220,-3){\line(1,0){30}}
        \put(255,-5){non $\alpha_3$-critical edges}
		\caption{\label{3}\small All vertices of $\{u'_1,\ldots,u'_p,v'_1,\ldots,v'_r\}$ that are static-excluded in $T$.}
		\end{picture}
	\end{center}
\end{figure}

Let $T_u$ (resp.~$T_v$) be the connected component of $T-uv$ containing $u$ (resp.~$v$). Let $\{u_1,\ldots,u_h\}$ (resp.~$\{v_1,\ldots,v_k\}$) be the vertices of $T_u$ (resp.~$T_v$) that are adjacent to $u$ (resp.~$v$) and are flexible in $T$, and let $\{e_1,\ldots,e_{h'}\}$ (resp.~$\{f_1,\ldots,f_{k'}\}$) be the $\alpha_3$-critical edges of $T$ that are incident to $u_i$ (resp.~$v_i$), as shown in Figure \ref{3}. Let $T'$ be the forest obtained by deleting all edges of $\{uu_1,\ldots,uu_h,vv_1,\ldots,vv_k\}$ from $T$.

By Corollary \ref{3.4}, all edges in $\{uv,e_1,\ldots,e_{h'},f_1,\ldots,f_{k'}\}$ are $\mu_3$-critical in $T$ and are covered by all maximum 3-matchings of $T$. Thus, none of edges in $\{uu_1,\ldots,uu_h,vv_1,\ldots,vv_k\}$ is covered by any maximum 3-matching of $T$, this means that $\mu_3(T')=\mu_3(T)$. By Theorem \ref{2.1}, $\alpha_3(T')+\mu_3(T')=|V(T')|=|V(T)|=\alpha_3(T)+\mu_3(T)$. Hence, $\alpha_3(T')=\alpha_3(T)$.

On the other hand, let $S$ be a maximum dissociation set of $T$ that does not contain $u$. Then $S\cup\{u\}$ is a dissociation set of $T'$, which implies that $\alpha_3(T')\geq \alpha_3(T)+1$, a contradiction.

\textbf{Case 2.} There is only one vertex in $\{u,v\}$ that has a neighbor in $\mathcal{A}_T$.

W.l.o.g, assume that $v$ has a neighbor $w$ in $\mathcal{A}_T$. The proofs are almost identical, the major change being the substitution of $T'=T-\{vw,uu_1,\ldots,uu_h,vv_1,\ldots,vv_k\}$ for $T'=T-\{uu_1,\ldots,uu_h,vv_1,\ldots,vv_k\}$. Since $w$ is a vertex in $\mathcal{A}_T$, the edge $vw$ is not $\alpha_3$-critical in $T$. Thus, $vw$ is not $\mu_3$-critical in $T$ and there exists a maximum $3$-matching of $T$ that does not cover $vw$. We have $\mu_3(T')=\mu_3(T)$ and $\alpha_3(T')=\alpha_3(T)$.

On the other hand, let $S$ be a maximum dissociation set of $T$ that does not contain $u$. It can be easily seen that $S\cup\{u\}$ is a dissociation set of $T'$. Thus, $\alpha_3(T')\geq \alpha_3(T)+1$ which leads to a contradiction.


(3) Let $v$ be any vertex in $\mathcal{A}_T$. We first show that $v$ is not saturated by all maximum 3-matchings of $T$. Suppose, to the contrary, that $v$ is saturated by all maximum 3-matchings of $T$, which means that $\mu_3(T-v)=\mu_3(T)-1$. Since both $T$ and $T-v$ are 3-K\"{o}nig-Egerv\'{a}ry graphs, $\alpha_3(T-v)+\mu_3(T-v)=n-1$ and $\alpha_3(T)+\mu_3(T)=n$. Thus
$\alpha_3(T)=\alpha_3(T-v)$. On the other hand, since $v$ is in $\mathcal{A}_T$, $\alpha_3(T-v)=\alpha_3(T)-1$, a contradiction.

Let $u_1u_2u_3$ be a critical $3$-path of $T$, and let $M$ be a maximum 3-matching of $T$ which does not saturate $v$. Suppose, to the contrary, that $v$ is adjacent to a vertex of $\{u_1,u_2,u_3\}$. W.l.o.g., suppose that $v$ is adjacent to $u_1$. Since $u_1u_2u_3$ is critical in $T$, it is a 3-path of $M$.
Then $M'=M\setminus\{u_1u_2u_3\}\cup\{vu_1u_2\}$ is also a maximum 3-matching of $T$ and the edge $u_2u_3$ is not covered by $M'$, which contradicts to the fact that $u_2u_3$ is $\mu_3$-critical in $T$ and it should be covered by any maximum 3-matching of $T$.

The proof is complete. \qed

The third structure theorem shows that every maximum dissociation set of a tree contains exactly one end-vertex of each insulated $\alpha_3$-critical edge and two vertices of each critical 3-path.

\begin{theo}\label{3.7}
Let $T$ be a tree and $\eta(T)$ the number of $\alpha_3$-critical edges of $T$. Then

(1) every maximum dissociation set contains exactly one end-vertex of each insulated $\alpha_3$-critical edge;

(2) every maximum dissociation set contains exactly two vertices of each critical $3$-path;

(3) ${\alpha _3}(T)=|\mathcal{A}_T|+\eta (T)$.

\end{theo}

\pf (1) By Theorem \ref{3.6}(2), every end-vertex of each insulated $\alpha_3$-critical edge has exactly one neighbor in $\mathcal{A}_T$. On the other hand, every maximum dissociation set contains all vertices of $\mathcal{A}_T$ and at least one end-vertex of each $\alpha_3$-critical edge. Thus, every maximum dissociation set contains exactly one end-vertex of each insulated $\alpha_3$-critical edge.

(2) Change the tree $T$ into a rooted tree by choosing any vertex as the root. Let
\begin{eqnarray*}
S=\mathcal{A}_T\cup \{u|uv~\text{is an}~\alpha_3\text{-critical edge of}~T ~\text{and~} \ell(u)>\ell(v)\}.
\end{eqnarray*}
If $uv$ is an edge of $T$, it is impossible that $\ell(u)=\ell(v)$. Thus, $S$ is well-defined. If $v$ is a vertex in the rooted tree other than the root, the parent of $v$ is the unique vertex $u$ such that there is a directed edge from $u$ to $v$. It's a simple fact that if $uv$ is an edge of $T$, then  $\ell(u)<\ell(v)$ if and only if $u$ is the parent of $v$ in the rooted tree. We divide our proof in three steps.

First, our task is to show that $S$ contains exactly two vertices of each critical 3-path. Clearly, $S$ contains at most two vertices of each critical 3-path. Let $u_1u_2u_3$ be a critical 3-path of $T$. Assume that $|S\cap\{u_1,u_2,u_3\}|=1$. According to the definition of $S$, $S\cap\{u_1,u_2,u_3\}=\{u_2\}$. Thus, $\ell(u_1)<\ell(u_2)$, $\ell(u_3)<\ell(u_2)$, and both $u_1$ and $u_3$ are the parents of $u_2$ in the rooted tree. This is impossible. Hence, $S$ contains exactly two vertices of each critical 3-path.

Next, we need to prove that $S$ is a dissociation set of $T$. To prove the assertion, we present the following claim.

\emph{Claim:} Let $v$ be a vertex of $\mathcal{F}_T$, and let $uv$ be an edge rather than an $\alpha_3$-critical edge of $T$. If $\ell(u)<\ell(v)$, then $S$ cannot contain the vertex $v$.

\textbf{Proof of Claim.} Let $T_v$ be the subtree rooted at $v$. If $w$ is a vertex in $T_v$ rather than the vertex $v$, then $\ell(w)>\ell(v)$. According to the definition of the set $S$, $S$ cannot contain the vertex $v$. \qed

By Theorem \ref{3.6} and Claim, if the induced subgraph $T[S]$ contains a 3-path, say $v_1v_2v_3$, then $\{v_1,v_3\}\subseteq\mathcal{F}_T$ and $v_2\in \mathcal{A}_T$. Moreover $\ell(v_1)<\ell(v_2)$, $\ell(v_3)<\ell(v_2)$ and both $v_1$ and $v_3$ are the parents of $v_2$ in the rooted tree. This is impossible. Thus we are led to the conclusion that $S$ is a dissociation set of $T$.

Finally, we prove the statement in (2). Let $F$ be a maximum dissociation set of $T$. Clearly, $\mathcal{A}_T\subseteq F\subseteq \mathcal{A}_T\cup\mathcal{F}_T$. By Theorem \ref{3.7}(1), $|F\cap \mathcal{F}_T|\leq \eta(T)$ and $|F|\leq|\mathcal{A}_T|+\eta(T)$. On the other hand, $S$ is a dissociation set and $|S|=|\mathcal{A}_T|+\eta(T)$. Thus, $S$ is a maximum dissociation set of $T$, and every maximum dissociation of $T$ contains exactly two vertices of each critical 3-path.

(3) Since the set $S$ defined in the proof of (2) is a maximum dissociation set of $T$ and $|S|=|\mathcal{A}_T|+\eta (T)$, ${\alpha _3}(T)=|\mathcal{A}_T|+\eta (T)$. \qed

The fourth structure theorem gives adjacency rules that determine the sets of static-included vertices and static-excluded vertices of a tree.

\begin{theo} \label{3.8}
Let $T$ be a tree and $u$ a vertex of $T$. If $u$ is in $\mathcal{N}_T$ and is adjacent to $p$ isolated vertices of $T[\mathcal{A}_T]$ and $q$ end-vertices of isolated edges of $T[\mathcal{A}_T]$, then $p+2q\geq 4$ or $p=3$. Thus, if $\mathcal{N}_T\neq\emptyset$, then $|\mathcal{A}_T|\geq3$.
\end{theo}

\pf Change the tree $T$ into a rooted tree by choosing the vertex $u$ as the root. Let
\begin{eqnarray*}
S=\mathcal{A}_T\cup \{v|vw~\text{is an}~\alpha_3\text{-critical edge of}~T ~\text{and~} \ell(v)>\ell(w)\}.
\end{eqnarray*}
According to the proof of Theorem \ref{3.7}(2), $S$ is a maximum dissociation set of $T$, and the vertex $u$ is not adjacent to any of the flexible vertices contained in $S$.

Next, we show that $p+2q\geq 4$ or $p=3$ and consider the following three cases.

\textbf{Case 1.} $q=0$

Let $N_T(u)\cap\mathcal{A}_T=\{u_1,\ldots,u_p\}$. Suppose, for a contradiction, that $p\leq2$. According to Claim in the proof of Theorem \ref{3.7}(2), each $u_i$ is not adjacent to any of the flexible vertices contained in $S$. Thus, $(S\setminus\{u_p\})\cup\{u\}$ is also a maximum dissociation set of $T$, a contradiction. Hence, in this case $p\geq3$.

\textbf{Case 2.} $q=1$

Let $w_1w_2$ be an edge of $T[\mathcal{A}_T]$ and $w_1$ a neighbor of $u$. Suppose, for a contradiction, that $p\leq 1$. Now, $(S\setminus\{w_1\})\cup\{u\}$ is also a maximum dissociation set of $T$. This leads to a contradiction. Hence, in this case $p\geq 2$.

\textbf{Case 3.} $q\geq2$

In this case, it is obvious that $p+2q\geq4$.

Consequently, we infer that $p+2q\geq 4$ or $p=3$. It follows that if $\mathcal{N}_T\neq \emptyset$, then $|\mathcal{A}_T|\geq3$. \qed

\section{The maximum number of maximum dissociation sets of a tree}

In this section, we use the structure theorems presented in Section 3 to find upper bounds on the number of maximum dissociation sets among all trees of order $n$.

\vspace{1mm}

\begin{lem} \label{4.1}
Let $T$ be a tree with $k$ flexible vertices. If $T$ has $x$ critical $3$-paths and $\frac{{k-3x}}{2}$ insulated $\alpha_3$-critical edges, then $T$ has at most ${3^x} \cdot {2^{\frac{{k-3x}}{2}}}$ maximum dissociation sets.
\end{lem}

\pf By Theorem \ref{3.6}, there is no critical $4$-path or critical $K_{1,3}$ in $T$. By Theorem \ref{3.7}, every maximum dissociation set contains exactly two vertices of each of the $x$ critical $3$-paths and exactly one end-vertex of each of the $\frac{{k-3x}}{2}$ insulated $\alpha_3$-critical edges. Thus, $T$ has at most ${3^x} \cdot {2^{\frac{{k - 3x}}{2}}}$ maximum dissociation sets. \qed


Let $f(x)={3^x}\cdot {2^{\frac{{k-3x}}{2}}}$. Since $f(x)$ is an increasing function of $x$,
we have

\begin{lem}\label{4.2}
Let $T$ be a tree with $k$ flexible vertices. Then

(1) if $k=3t$, $T$ has at most $3^t$ maximum dissociation sets, and the upper bound is achieved only if $T$ contains $t$ critical 3-paths,

(2) if $k=3t+1$, $T$ has at most $3^{t-1}\cdot 2^2$ maximum dissociation sets, and the upper bound is achieved only if $T$ contains $t-1$ critical 3-paths and two insulated $\alpha_3$-critical edges,

(3) if $k=3t+2$, $T$ has at most $3^t\cdot 2$ maximum dissociation sets, and the upper bound is achieved only if $T$ contains $t$ critical 3-paths and one insulated $\alpha_3$-critical edge,

(4) the larger $k$, the larger upper bound of the number of maximum dissociation sets of $T$.
\end{lem}
%

\begin{figure}[h]
	\begin{center}
		\begin{picture}(400,40)
		\put(220,40){\line(1,0){10}}
		\put(220,40){\line(0,-1){30}}
		\put(230,40){\line(0,-1){30}}
		\put(220,10){\line(1,0){10}}
     	\put(270,30){\line(1,0){30}}
		\put(270,30){\line(0,-1){10}}
		\put(300,30){\line(0,-1){10}}
		\put(270,20){\line(1,0){30}}
        \put(245,20){\emph{or}}
		\put(90,-3){A critical 3-path}
       \multiput(100,20)(30,0){3}{\circle*{3.5}}
        \linethickness{0.5mm}
		\put(100,20){\line(1,0){60}}
		\caption{\label{4}\small A critical 3-path and its symbolic representations.}
		\end{picture}
	\end{center}
\end{figure}
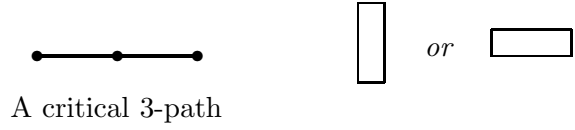

Let $S^*_{T_1,\ldots,T_r}$ be the tree consisting of $r$ induced subtrees $T_1, \ldots$, $T_r$ with a
common leaf. For simplicity, we use a rectangle to represent a critical 3-path (see Figure \ref{4}). A vertex $v$ is said to be adjacent to a critical 3-path $P$ if $v$ is adjacent to a vertex of $P$. Two critical 3-paths are adjacent if they are connected by an edge.

\begin{theo}\label{4.3}
Let $m$ be a positive integer and $T$ a tree on $3m+1$ vertices. Then $T$ has at most ${3^{m-1}} + 1$ maximum dissociation sets. This bound is best possible. When $m\geq2$, the bound is achieved only on the families of trees $S^*_{P_3,P_2,T_1,\ldots,T_{m-1}}$, where $T_i\cong P_4\text{\ or\ } K_{1,3}$ ($1\leq i\leq m-1$) (see Figure \ref{5}).
\end{theo}
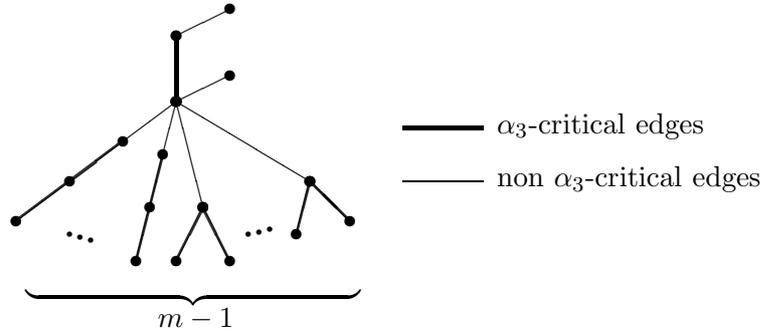
\begin{figure}[h]
	\begin{center}
		\begin{picture}(400,120)
		\multiput(170,80)(-5,-20){4}{\circle*{3.5}}
		\multiput(170,80)(-20,-15){4}{\circle*{3.5}}
		\multiput(170,80)(10,-40){2}{\circle*{3.5}}
		\multiput(170,80)(50,-30){2}{\circle*{3.5}}
		\multiput(180,40)(-10,-20){2}{\circle*{3.5}}
		\multiput(180,40)(10,-20){2}{\circle*{3.5}}
		\multiput(220,50)(-5,-20){2}{\circle*{3.5}}
		\multiput(220,50)(15,-15){2}{\circle*{3.5}}
		\multiput(170,80)(0,25){2}{\circle*{3.5}}
		\multiput(170,80)(20,10){2}{\circle*{3.5}}
		\multiput(170,105)(20,10){2}{\circle*{3.5}}
		\put(170,80){\line(-1,-4){15}}
		\put(170,80){\line(-4,-3){60}}
		\put(170,80){\line(5,-3){50}}
		\put(180,40){\line(-1,4){10}}
		\put(170,80){\line(2,1){20}}
		\put(170,105){\line(2,1){20}}
		\put(165.2,59.8){\line(-1,-4){10}}
		\put(164.8,60.2){\line(-1,-4){10}}
		\put(150.2,64.8){\line(-4,-3){40}}
		\put(149.8,65.2){\line(-4,-3){40}}
		\put(180,40){\line(-1,-2){10}}
		\put(180.2,39.8){\line(-1,-2){10}}
		\put(179.8,40.2){\line(-1,-2){10}}
		\put(180,40){\line(1,-2){10}}
		\put(180.2,40.2){\line(1,-2){10}}
		\put(179.8,39.8){\line(1,-2){10}}
		\put(220,50){\line(-1,-4){5}}
		\put(220.2,49.8){\line(-1,-4){5}}
		\put(219.8,50.2){\line(-1,-4){5}}
		\put(220,50){\line(1,-1){15}}
		\put(220.2,50.2){\line(1,-1){15}}
		\put(219.8,49.8){\line(1,-1){15}}
		\put(110,10){\ $\underbrace {\kern 125pt}$}
		\put(163,-5){$m-1$}
		\linethickness{0.5mm}
		\put(255,70){\line(1,0){30}}
		\put(170,80){\line(0,1){25}}
		\thinlines
		\put(255,50){\line(1,0){30}}
		\put(290,68){$\alpha_3$-critical edges}
		\put(290,48){non $\alpha_3$-critical edges}
		\multiput(197,30)(4,1){3}{\circle*{2}}
		\multiput(130,30)(4,-1){3}{\circle*{2}}
		\caption{\label{5}\small Trees on $3m+1$ vertices with the most maximum dissociation sets. }
		\end{picture}
	\end{center}
\end{figure}

\pf We proceed to prove this theorem and distinguish the following two cases.

\textbf{Case 1.} $\mathcal{N}_T\neq\emptyset$.

By Theorem \ref{3.8}, $|\mathcal{A}_T|\geq3$. Thus, there are at most $3m-3$ flexible vertices in $T$. By Lemma \ref{4.2}, in this case $T$ has at most $3^{m-1}$ maximum dissociation sets.

\textbf{Case 2.} $\mathcal{N}_T=\emptyset$.

We first show that $|\mathcal{A}_T|\geq 2$. Suppose, for a contradiction, that $|\mathcal{A}_T|\leq1$. By Theorem \ref{3.6}, there is no insulated $\alpha_3$-critical edge in $T$, and $\mathcal{A}_T=\emptyset$. Now, each vertex in $T$ is a vertex of a critical $3$-path, which contradicts the fact that the number of vertices of $T$ is $3m+1$. Thus, we have proved that $|\mathcal{A}_T|\geq 2$.

Next, we consider the following three subcases.

\textbf{Subcase 2.1.} $|\mathcal{A}_T|=2$.

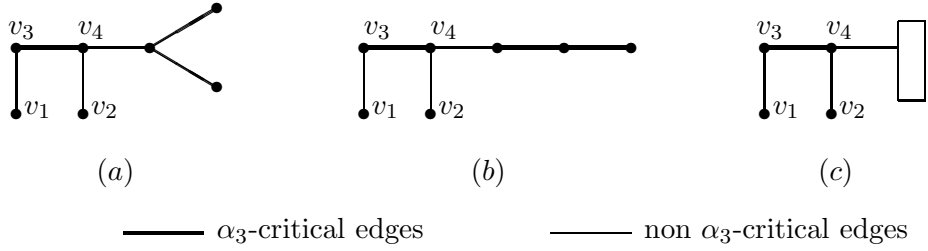
\begin{figure}[h]
	\begin{center}
		\begin{picture}(400,100)
		\multiput(150,75)(25,0){5}{\circle*{3.5}}
		\multiput(150,50)(25,0){2}{\circle*{3.5}}
		\linethickness{0.4mm}
		\put(150,75){\line(1,0){25}}
		\put(200,75){\line(1,0){50}}
		\put(20,75){\line(1,0){25}}
		\thinlines
		\put(150,75){\line(0,-1){25}}
		\put(175,75){\line(1,0){25}}
		\put(175,75){\line(0,-1){25}}
		\put(20,75){\line(0,-1){25}}
		\put(45,75){\line(0,-1){25}}
		\put(45,75){\line(1,0){25}}
		\put(70,75){\line(5,3){25}}
		\put(70.2,74.8){\line(5,3){25}}
		\put(69.8,75.2){\line(5,3){25}}
		\put(70,75){\line(5,-3){25}}
		\put(70.2,75.2){\line(5,-3){25}}
		\put(69.8,74.8){\line(5,-3){25}}
		\put(150,80){$v_3$}
		\put(175,80){$v_4$}
		\put(153,50){$v_1$}
		\put(178,50){$v_2$}
		\put(23,50){$v_1$}
		\put(48,50){$v_2$}
		\put(17,80){$v_3$}
		\put(43,80){$v_4$}
		\put(190,25){($b$)}
		\put(50,25){($a$)}
		\put(320,25){($c$)}
		\multiput(20,75)(25,0){3}{\circle*{3.5}}
		\multiput(20,50)(25,0){2}{\circle*{3.5}}
		\multiput(95,90)(0,-30){2}{\circle*{3.5}}
		\linethickness{0.4mm}
		\put(300,75){\line(1,0){25}}
		\thinlines
		\put(300,75){\line(0,-1){25}}
		\put(325,75){\line(0,-1){25}}
		\put(325,75){\line(1,0){25}}
		\put(350,85){\line(1,0){10}}
		\put(350,85){\line(0,-1){30}}
		\put(360,85){\line(0,-1){30}}
		\put(350,55){\line(1,0){10}}
		\put(303,50){$v_1$}
		\put(328,50){$v_2$}
		\put(297,80){$v_3$}
		\put(323,80){$v_4$}
		\multiput(300,75)(25,0){2}{\circle*{3.5}}
		\multiput(300,50)(25,0){2}{\circle*{3.5}}
		\linethickness{0.4mm}
		\put(60,5){\line(1,0){30}}
        \put(95,3){$\alpha_3$-critical edges}
		\thinlines
		\put(220,5){\line(1,0){30}}
		\put(255,3){non $\alpha_3$-critical edges}
		\caption{\label{6}\small Subgraphs used in Subcase 2.1, where (c) is a symbolic representation of (a) and (b).}
		\end{picture}
	\end{center}
\end{figure}

By Theorem \ref{3.6}, $T$ has one insulated $\alpha_3$-critical edge and $m-1$ critical 3-paths. Let $\mathcal{A}_T=\{v_1,v_2\}$, and let $v_3v_4$ be the insulated $\alpha_3$-critical edge of $T$. W.l.o.g, suppose that the vertex $v_4$ is adjacent to a critical 3-path, as shown in Figure \ref{6}.

When $m=2$, any tree in the family pictured in Figure \ref{5} has $3^{m-1}+1$ maximum dissociation sets.

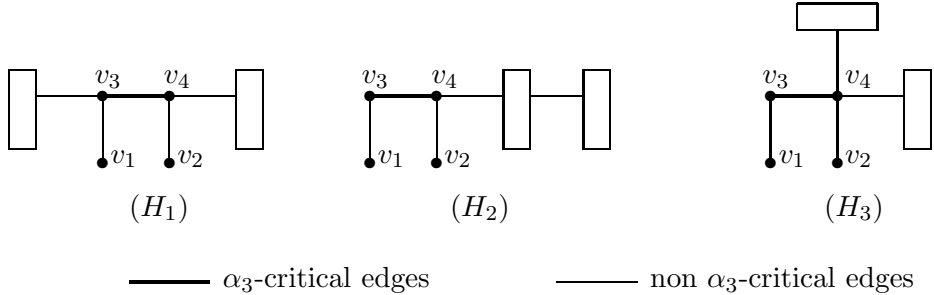
\begin{figure}[h]
	\begin{center}
		\begin{picture}(400,110)
		\linethickness{0.4mm}
		\put(50,75){\line(1,0){25}}
		\thinlines
		\put(50,75){\line(0,-1){25}}
		\put(75,75){\line(0,-1){25}}
		\put(75,75){\line(1,0){25}}
		\put(100,85){\line(1,0){10}}
		\put(100,85){\line(0,-1){30}}
		\put(110,85){\line(0,-1){30}}
		\put(100,55){\line(1,0){10}}
		\put(50,75){\line(-1,0){25}}
		\put(15,85){\line(1,0){10}}
		\put(15,85){\line(0,-1){30}}
		\put(25,85){\line(0,-1){30}}
		\put(15,55){\line(1,0){10}}
		\put(53,50){$v_1$}
		\put(78,50){$v_2$}
		\put(47,80){$v_3$}
		\put(73,80){$v_4$}
		\put(60,30){($H_1$)}
		\multiput(50,75)(25,0){2}{\circle*{3.5}}
		\multiput(50,50)(25,0){2}{\circle*{3.5}}
		\linethickness{0.4mm}
		\put(60,5){\line(1,0){30}}
        \put(95,3){$\alpha_3$-critical edges}
		\thinlines
		\put(220,5){\line(1,0){30}}
		\put(255,3){non $\alpha_3$-critical edges}
		\linethickness{0.4mm}
		\put(150,75){\line(1,0){25}}
		\thinlines
		\put(150,75){\line(0,-1){25}}
		\put(175,75){\line(0,-1){25}}
		\put(175,75){\line(1,0){25}}
		\put(200,85){\line(1,0){10}}
		\put(200,85){\line(0,-1){30}}
		\put(210,85){\line(0,-1){30}}
		\put(200,55){\line(1,0){10}}
		\put(210,75){\line(1,0){20}}
		\put(230,85){\line(1,0){10}}
		\put(230,85){\line(0,-1){30}}
		\put(240,85){\line(0,-1){30}}
		\put(230,55){\line(1,0){10}}
		\put(153,50){$v_1$}
		\put(178,50){$v_2$}
		\put(147,80){$v_3$}
		\put(173,80){$v_4$}
		\put(180,30){($H_2$)}
		\multiput(150,75)(25,0){2}{\circle*{3.5}}
		\multiput(150,50)(25,0){2}{\circle*{3.5}}
		\linethickness{0.4mm}
		\put(300,75){\line(1,0){25}}
		\thinlines
		\put(300,75){\line(0,-1){25}}
		\put(325,75){\line(0,-1){25}}
		\put(325,75){\line(1,0){25}}
		\put(350,85){\line(1,0){10}}
		\put(350,85){\line(0,-1){30}}
		\put(360,85){\line(0,-1){30}}
		\put(350,55){\line(1,0){10}}
		\put(325,75){\line(0,1){25}}
		\put(340,100){\line(0,1){10}}
		\put(310,100){\line(0,1){10}}
		\put(310,100){\line(1,0){30}}
		\put(310,110){\line(1,0){30}}
		\put(303,50){$v_1$}
		\put(328,50){$v_2$}
		\put(297,80){$v_3$}
		\put(328,80){$v_4$}
		\put(320,30){($H_3$)}
		\multiput(300,75)(25,0){2}{\circle*{3.5}}
		\multiput(300,50)(25,0){2}{\circle*{3.5}}
		\caption{\label{7}\small Three families $H_1$,  $H_2$, and $H_3$.}
		\end{picture}
	\end{center}
\end{figure}

When $m\geq 3$, $T$ must contain a subgraph that is isomorphic to a tree in the family $H_1$, or to a tree in the family $H_2$, or to a tree in the family $H_3$. The three families of trees $H_1$, $H_2$ and $H_3$ are shown in Figure \ref{7}. A tree in the family $H_1$ has exactly 6 maximum dissociation sets. Thus, if $T$ contains a subgraph that is isomorphic to a tree in the family $H_1$, then $T$ has at most $6\cdot 3^{m-3}$ maximum dissociation sets. A tree in the family $H_2$ has at most 9 maximum dissociation sets. Thus if $T$ contains a subgraph that is isomorphic to a tree in the family $H_2$, then $T$ has at most $9\cdot 3^{m-3}$ maximum dissociation sets.

It follows that in this subcase if $T$ has $3^{m-1}+1$ maximum dissociation sets, then all critical 3-paths of $T$ are adjacent to a common vertex. On the other hand, consider a tree $T$ in which all critical 3-paths are adjacent to a common vertex, i.e., a tree in the families pictured in Figure \ref{5}. By simple calculation, there are $3^{m-1}+1$ maximum dissociation sets in $T$. Thus, the upper bound can be achieved by these families of trees.

\textbf{Subcase 2.2.} $|\mathcal{A}_T|=3$.

\begin{figure}[h]
	\begin{center}
		\begin{picture}(400,30)
		\multiput(100,30)(30,0){2}{\circle*{3.5}}
		\multiput(90,5)(50,0){3}{\circle*{3.5}}
		\multiput(150,30)(30,0){2}{\circle*{3.5}}
		\linethickness{0.4mm}
		\put(100,30){\line(1,0){30}}
		\put(150,30){\line(1,0){30}}
		\thinlines
		\put(100,30){\line(-2,-5){10}}
		\put(130,30){\line(2,-5){10}}	
		\put(150,30){\line(-2,-5){10}}
		\put(180,30){\line(2,-5){10}}
		\linethickness{0.4mm}
		\put(225,25){\line(1,0){30}}
		\put(262,25){$\alpha_3$-critical edges}
		\thinlines
		\put(225,3){\line(1,0){30}}
		\put(262,3){non $\alpha_3$-critical edges}
		\caption{\label{8}\small A subgraph $H$ that is isomorphic to $P_7$.}
		\end{picture}
	\end{center}
\end{figure}
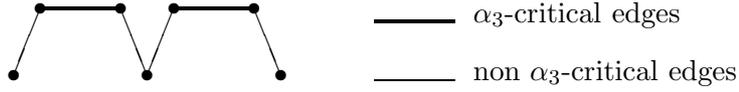

By Theorem \ref{3.6}, $T$ has two insulated $\alpha_3$-critical edges and $m-2$ critical 3-paths. Moreover, $T$ contains a subgraph $H$ that is isomorphic to $P_7$. See Figure \ref{8}. Since $H$ has three maximum dissociation sets, in this subcase $T$ has at most $3\cdot 3^{m-2}$ maximum dissociation sets.

\textbf{Subcase 2.3.} $|\mathcal{A}_T|\geq4$.

In this subcase, there are at most $3m-3$ flexible vertices in $T$. By Lemma \ref{4.2}, $T$ has at most $3^{m-1}$ maximum dissociation sets.

We can now derive the final conclusion. A tree on $3m+1$ vertices has at most $3^{m-1}+1$ maximum dissociation sets. When $m\geq2$, the upper bound is achieved only in Subcase 2.1 and only on the families of trees pictured in Figure \ref{5}. \qed

\begin{figure}[h]
	\begin{center}
		\begin{picture}(400,35)
		\linethickness{0.4mm}
		\put(210,25){\line(1,0){30}}
		\put(245,23){$\alpha_3$-critical edges}
		\thinlines
		\put(210,1){\line(1,0){30}}
		\put(245,-1){non $\alpha_3$-critical edges}
		\linethickness{0.4mm}
		\put(100,25){\line(1,0){25}}
		\put(150,25){\line(1,0){25}}
		\thinlines
		\put(100,25){\line(0,-1){25}}
		\put(125,25){\line(0,-1){25}}
		\put(150,25){\line(0,-1){25}}
		\put(175,25){\line(0,-1){25}}
		\put(125,25){\line(1,0){25}}
		\put(103,0){$v_1$}
		\put(128,0){$v_2$}
		\put(97,30){$u_1$}
		\put(123,30){$u_2$}
		\put(153,0){$v_3$}
		\put(178,0){$v_4$}
		\put(147,30){$u_3$}
		\put(173,30){$u_4$}
		\multiput(100,25)(25,0){4}{\circle*{3.5}}
		\multiput(100,0)(25,0){4}{\circle*{3.5}}
		\caption{\label{9}\small A special tree $LT_8$ on 8 vertices.}
		\end{picture}
	\end{center}
\end{figure}
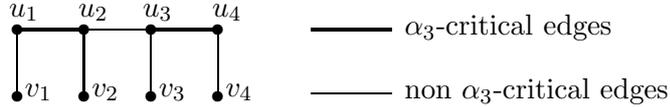

We now handle the case when the number of vertices of $T$ is $3m+2$. Define a special tree $LT_8$ on 8 vertices to contain two insulated $\alpha_3$-critical edges $u_1u_2$ and $u_3u_4$, and five non $\alpha_3$-critical edges $u_2u_3$, and $u_iv_i$ for $i=1$ to 4. See Figure \ref{9}.

\begin{theo}\label{4.4}
Let $m$ be a positive integer and $T$ a tree on $3m+2$ vertices. Then $T$ has at most $3^{m-1}$ maximum dissociation sets. This bound is best possible. When $m=2$, the bound is achieved only on the families of trees $LT_8$, $S^*_{P_3,P_3,T_1}$, $S^*_{P_2,P_2,P_2,P_2,T_1}$, and $S^*_{P_3,P_2,P_2,T_1}$, where $T_1\cong P_4 \text{\ or\ }K_{1,3}$. When $m\neq 2$, the bound is achieved only on the families of trees $S^*_{P_3,P_3,T_1,\ldots,T_{m-1}}$,  $S^*_{P_2,P_2,P_2,P_2,T_1,\ldots,T_{m-1}}$ and $S^*_{P_3,P_2,P_2,T_1,\ldots,T_{m-1}}$ where $T_i\cong P_4 \text{\ or\ } K_{1,3}$ ($1\leq i\leq m-1$) (see Figure \ref{10}).
\end{theo}
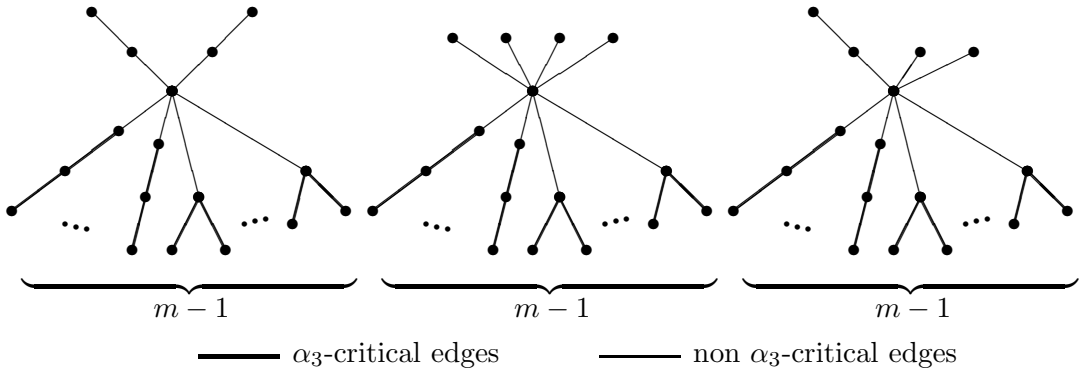
\begin{figure}[h]
	\begin{center}
		\begin{picture}(400,120)
		\multiput(195,90)(-5,-20){4}{\circle*{3.5}}
		\multiput(195,90)(-20,-15){4}{\circle*{3.5}}
		\multiput(195,90)(10,-40){2}{\circle*{3.5}}
		\multiput(195,90)(50,-30){2}{\circle*{3.5}}
		\multiput(205,50)(-10,-20){2}{\circle*{3.5}}
		\multiput(205,50)(10,-20){2}{\circle*{3.5}}
		\multiput(245,60)(-5,-20){2}{\circle*{3.5}}
		\multiput(245,60)(15,-15){2}{\circle*{3.5}}
		\multiput(165,110)(20,0){4}{\circle*{3.5}}
		\put(195,90){\line(-1,-4){15}}
		\put(195,90){\line(-4,-3){60}}
		\put(195,90){\line(5,-3){50}}
		\put(205,50){\line(-1,4){10}}
		\put(190.2,69.8){\line(-1,-4){10}}
		\put(189.8,70.2){\line(-1,-4){10}}
		\put(175.2,74.8){\line(-4,-3){40}}
		\put(174.8,75.2){\line(-4,-3){40}}
		\put(205,50){\line(-1,-2){10}}
		\put(205.2,49.8){\line(-1,-2){10}}
		\put(204.8,50.2){\line(-1,-2){10}}
		\put(205,50){\line(1,-2){10}}
		\put(205.2,50.2){\line(1,-2){10}}
		\put(204.8,49.8){\line(1,-2){10}}
		\put(245,60){\line(-1,-4){5}}
		\put(245.2,59.8){\line(-1,-4){5}}
		\put(244.8,60.2){\line(-1,-4){5}}
		\put(245,60){\line(1,-1){15}}
		\put(245.2,60.2){\line(1,-1){15}}
		\put(244.8,59.8){\line(1,-1){15}}
		\put(135,20){\ $\underbrace {\kern 125pt}$}
		\put(188,5){$m-1$}
		\multiput(222,40)(4,1){3}{\circle*{2}}
		\multiput(155,40)(4,-1){3}{\circle*{2}}
		\put(195,90){\line(-1,2){10}}
		\put(195,90){\line(-3,2){30}}
		\put(195,90){\line(1,2){10}}
		\put(195,90){\line(3,2){30}}
		\multiput(60,90)(-5,-20){4}{\circle*{3.5}}
		\multiput(60,90)(-20,-15){4}{\circle*{3.5}}
		\multiput(60,90)(10,-40){2}{\circle*{3.5}}
		\multiput(60,90)(50,-30){2}{\circle*{3.5}}
		\multiput(70,50)(-10,-20){2}{\circle*{3.5}}
		\multiput(70,50)(10,-20){2}{\circle*{3.5}}
		\multiput(110,60)(-5,-20){2}{\circle*{3.5}}
		\multiput(110,60)(15,-15){2}{\circle*{3.5}}
		\multiput(60,90)(-15,15){3}{\circle*{3.5}}
		\multiput(60,90)(15,15){3}{\circle*{3.5}}
		\put(60,90){\line(-1,-4){15}}
		\put(60,90){\line(-4,-3){60}}
		\put(60,90){\line(5,-3){50}}
		\put(70,50){\line(-1,4){10}}
		\put(55.2,69.8){\line(-1,-4){10}}
		\put(54.8,70.2){\line(-1,-4){10}}
		\put(40.2,74.8){\line(-4,-3){40}}
		\put(39.8,75.2){\line(-4,-3){40}}
		\put(70,50){\line(-1,-2){10}}
		\put(70.2,49.8){\line(-1,-2){10}}
		\put(69.8,50.2){\line(-1,-2){10}}
		\put(70,50){\line(1,-2){10}}
		\put(70.2,50.2){\line(1,-2){10}}
		\put(69.8,49.8){\line(1,-2){10}}
		\put(110,60){\line(-1,-4){5}}
		\put(110.2,59.8){\line(-1,-4){5}}
		\put(109.8,60.2){\line(-1,-4){5}}
		\put(110,60){\line(1,-1){15}}
		\put(110.2,60.2){\line(1,-1){15}}
		\put(109.8,59.8){\line(1,-1){15}}
		\put(0,20){\ $\underbrace {\kern 125pt}$}
		\put(53,5){$m-1$}
		\multiput(87,40)(4,1){3}{\circle*{2}}
		\multiput(20,40)(4,-1){3}{\circle*{2}}
		\put(60,90){\line(-1,1){30}}
		\put(60,90){\line(1,1){30}}
		\multiput(330,90)(-5,-20){4}{\circle*{3.5}}
		\multiput(330,90)(-20,-15){4}{\circle*{3.5}}
		\multiput(330,90)(10,-40){2}{\circle*{3.5}}
		\multiput(330,90)(50,-30){2}{\circle*{3.5}}
		\multiput(340,50)(-10,-20){2}{\circle*{3.5}}
		\multiput(340,50)(10,-20){2}{\circle*{3.5}}
		\multiput(380,60)(-5,-20){2}{\circle*{3.5}}
		\multiput(380,60)(15,-15){2}{\circle*{3.5}}
		\multiput(330,90)(-15,15){3}{\circle*{3.5}}
		\multiput(340,105)(20,0){2}{\circle*{3.5}}
		\put(330,90){\line(-1,1){30}}
		\put(330,90){\line(2,3){10}}
		\put(330,90){\line(2,1){30}}
		\put(330,90){\line(-1,-4){15}}
		\put(330,90){\line(-4,-3){60}}
		\put(330,90){\line(5,-3){50}}
		\put(340,50){\line(-1,4){10}}
		\put(325.2,69.8){\line(-1,-4){10}}
		\put(324.8,70.2){\line(-1,-4){10}}
		\put(310.2,74.8){\line(-4,-3){40}}
		\put(309.8,75.2){\line(-4,-3){40}}
		\put(340,50){\line(-1,-2){10}}
		\put(340.2,49.8){\line(-1,-2){10}}
		\put(339.8,50.2){\line(-1,-2){10}}
		\put(340,50){\line(1,-2){10}}
		\put(340.2,50.2){\line(1,-2){10}}
		\put(339.8,49.8){\line(1,-2){10}}
		\put(380,60){\line(-1,-4){5}}
		\put(380.2,59.8){\line(-1,-4){5}}
		\put(379.8,60.2){\line(-1,-4){5}}
		\put(380,60){\line(1,-1){15}}
		\put(380.2,60.2){\line(1,-1){15}}
		\put(379.8,59.8){\line(1,-1){15}}
		\put(270,20){\ $\underbrace {\kern 125pt}$}
		\put(323,5){$m-1$}
		\multiput(357,40)(4,1){3}{\circle*{2}}
		\multiput(290,40)(4,-1){3}{\circle*{2}}
		\linethickness{0.4mm}
		\put(70,-10){\line(1,0){30}}
		\thinlines
		\put(220,-10){\line(1,0){30}}
		\put(105,-12){$\alpha_3$-critical edges}
		\put(255,-12){non $\alpha_3$-critical edges}
		\caption{\label{10}\small Trees on $3m+2$ vertices with the most maximum dissociation sets.}
		\end{picture}
	\end{center}
\end{figure}

\pf We proceed to prove this theorem and distinguish the following two cases.

\textbf{Case 1.} $\mathcal{N}_T\neq\emptyset$.

It follows from Theorem \ref{3.8} that $|\mathcal{A}_T|\geq3$.

\textbf{Subcase 1.1.} $|\mathcal{A}_T|=3$.

Let $u$ be a vertex in $\mathcal{N}_T$, and let $\mathcal{A}_T=\{v_1,v_2,v_3\}$. By Theorem \ref{3.8}, all vertices of $\mathcal{A}_T$ are isolated vertices of $T[\mathcal{A}_T]$ and $T[\{u,v_1,v_2,v_3\}]$ is a star of four vertices.
We first show that there is no insulated $\alpha_3$-critical edge in $T$. Suppose, for a contradiction, that $e$ is an insulated $\alpha_3$-critical edge. By Theorem \ref{3.6}, every end-vertex of $e$ has a neighbor in $\{v_1,v_2,v_3\}$, which leads to a cycle in $T$. This is impossible. It follows that every vertex in $V(T)\setminus\{u,v_1,v_2,v_3\}$ is a vertex of a critical 3-path. This also contradicts the fact that the number of vertices of $V(T)\setminus\{u,v_1,v_2,v_3\}$ is $3m-2$.

Consequently,  this subcase is impossible to happen.

\textbf{Subcase 1.2.} $|\mathcal{A}_T|=4$.

In this subcase, $T$ has at most $3m-3$ flexible vertices. By Lemma \ref{4.2}, $T$ has at most $3^{m-1}$ maximum dissociation sets. And the upper bound is achieved only if $T$ contains one static-excluded vertex and $m-1$ critical 3-paths. By Theorem \ref{3.6}, all critical 3-paths are adjacent to the static-excluded vertex. Thus, the upper bound is achieved only if $T$ is a tree in the families pictured in Figure \ref{10}.
On the other hand, it can easily be seen that there are $3^{m-1}$ maximum dissociation sets in a tree in the families pictured in Figure \ref{10}. Thus, the upper bound can be achieved by these families of trees.

\textbf{Subcase 1.3.} $|\mathcal{A}_T|\geq5$.

It follows that there are at most $3m-4$ flexible vertices in $T$. By Lemma \ref{4.2}, in this subcase $T$ has at most $2\cdot 3^{m-2}$ maximum dissociation sets.

\textbf{Case 2.} $\mathcal{N}_T=\emptyset$.

We first show that $|\mathcal{A}_T|\geq 2$. Suppose, for a contradiction, that $|\mathcal{A}_T|\leq1$. By Theorem \ref{3.6}(2), there is no insulated $\alpha_3$-critical edge in $T$. Thus, $\mathcal{A}_T=\emptyset$. Now, every vertex in $T$ is a vertex of a critical 3-path. This contradicts the fact that the number of vertices of $T$ is $3m+2$. Hence, $|\mathcal{A}_T|\geq 2$.

\textbf{Subcase 2.1.} $|\mathcal{A}_T|=2$.

Let $\mathcal{A}_T=\{u_1,u_2\}$. By Theorem \ref{3.6}, in this subcase there is exactly one insulated $\alpha_3$-critical edge, say $v_1v_2$, in $T$. Let
$U=\{u_1,u_2,v_1,v_2\}$. Every vertex in $V(T)\setminus U$ is a vertex of a critical 3-path. This leads to a contradiction because the number of vertices of $V(T)\setminus U$ is $3m-2$. It follows that this subcase is impossible to happen.

\textbf{Subcase 2.2.} $|\mathcal{A}_T|=3$.

Let $\mathcal{A}_T=\{u_1,u_2,u_3\}$. By Theorem \ref{3.6}, in this subcase there are exactly two insulated $\alpha_3$-critical edges, say $v_1v_2$ and $v_3v_4$, in $T$. Let $U=\{u_1,u_2,u_3,v_1,v_2,v_3,v_4\}$. Every vertex in $V(T)\setminus U$ is a vertex of a critical 3-path. This leads to a contradiction because the number of vertices of $V(T)\setminus U$ is $3m-5$. It follows that this subcase is also impossible to happen.

\textbf{Subcase 2.3.} $|\mathcal{A}_T|=4$.

When $m=1$, this subcase is impossible to happen. When $m=2$, $T$ has 8 vertices and is isomorphic to $LT_8$. On the other hand, $LT_8$ has $3(=3^{m-1})$ maximum dissociation sets. Thus, when $m=2$, the upper bound can be achieved by the special tree $LT_8$.

\begin{figure}[h]
	\begin{center}
		\begin{picture}(200,40)
		\linethickness{0.4mm}
		\put(0,25){\line(1,0){25}}
		\thinlines
		\put(0,25){\line(0,-1){25}}
		\put(25,25){\line(0,-1){25}}
		\put(25,25){\line(1,0){25}}
		\put(50,35){\line(1,0){10}}
		\put(50,35){\line(0,-1){30}}
		\put(60,35){\line(0,-1){30}}
		\put(50,5){\line(1,0){10}}
		\put(3,0){$v_1$}
		\put(28,0){$v_2$}
		\put(-3,30){$v_3$}
		\put(23,30){$v_4$}
		\multiput(0,25)(25,0){2}{\circle*{3.5}}
		\multiput(0,0)(25,0){2}{\circle*{3.5}}

		\linethickness{0.4mm}
		\put(90,30){\line(1,0){30}}
        \put(130,28){$\alpha_3$-critical edges}
		\thinlines
		\put(90,6){\line(1,0){30}}
		\put(130,4){non $\alpha_3$-critical edges}
        \put(15,-15){$H$}
		\caption{\label{11}\small The family $H$.}
		\end{picture}
	\end{center}
\end{figure}
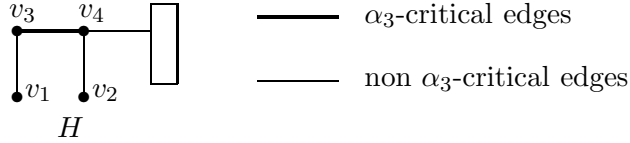

When $m>2$, there are at least two insulated $\alpha_3$-critical edges and at most $m-2$ critical 3-paths. Moreover, $T$ must contain a subgraph that is isomorphic to a tree in the family $H$ pictured in Figure \ref{11}. Since any tree in the family $H$ has at most 4 maximum dissociation sets, $T$ has at most $4\cdot 2^{2-1}\cdot 3^{(m-2)-1}$ maximum dissociation sets. It follows that in this subcase $T$ has at most $8\cdot 3^{m-3}$ maximum dissociation sets.

\textbf{Subcase 2.4.} $|\mathcal{A}_T|\geq 5$.

In this subcase, there are at most $3m-3$ flexible vertices and at least three insulated $\alpha_3$-critical edges in $T$. By Lemma \ref{4.2}, $T$ has at most $2^3\cdot 3^{m-3}$ maximum dissociation sets.

Consequently, we infer that a tree on $3m+2$ vertices has at most $3^{m-1}$ maximum dissociation sets. When $m=2$, the upper bound is achieved only in Subcase 1.2 and 2.3 and only on $LT_8$ and the families of trees pictured in Figure \ref{10}. When $m\neq 2$, the upper bound is achieved only in Subcase 1.2 and only on the families of trees pictured in Figure \ref{10}. \qed

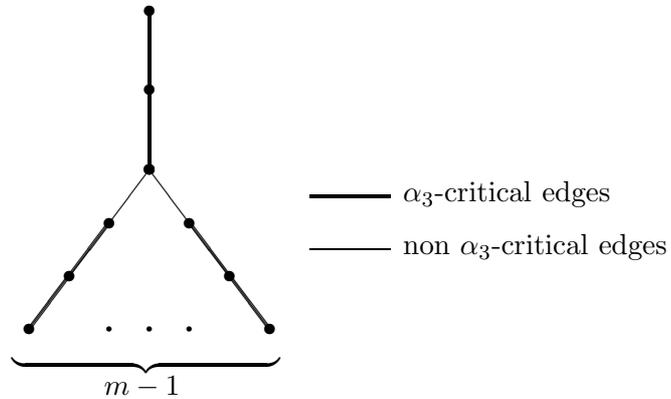
\begin{figure}[h]
	\begin{center}
		\begin{picture}(400,130)
		\multiput(200,80)(15,-20){4}{\circle*{3.5}}
		\multiput(200,80)(-15,-20){4}{\circle*{3.5}}
		\multiput(200,80)(0,30){3}{\circle*{3.5}}
		\put(200,80){\line(3,-4){45}}
		\put(200,80){\line(-3,-4){45}}
		\put(200.4,80){\line(0,1){60}}
		\put(200,80){\line(0,1){60}}
		\put(199.6,80){\line(0,1){60}}
		\put(215.3,60.3){\line(3,-4){30}}
		\put(214.7,59.7){\line(3,-4){30}}
		\put(185.3,59.7){\line(-3,-4){30}}
		\put(184.7,60.3){\line(-3,-4){30}}
		\put(145,10){\ $\underbrace {\kern 100pt}$}
		\put(183,-5){$m-1$}
		\linethickness{0.4mm}
		\put(260,70){\line(1,0){30}}
		\thinlines
		\put(260,50){\line(1,0){30}}
		\put(295,68){$\alpha_3$-critical edges}
		\put(295,48){non $\alpha_3$-critical edges}
		\multiput(185,20)(15,0){3}{\circle*{2}}
		\caption{\label{12} \small Trees on $3m$ vertices with the most maximum dissociation sets.}
		\end{picture}
	\end{center}
\end{figure}

\begin{theo} \label{4.5}
Let $m$ be a positive integer and $T$ a tree on $3m$ vertices. Then $T$ has at most ${3^{m-1}}+m+1$ maximum dissociation sets. This bound is best possible. The bound is achieved only on the family of trees $S^*_{P_3,T_1,\ldots,T_{m-1}}$ where $T_i\cong P_4, (1\leq i\leq m-1)$ (see Figure \ref{12}).
\end{theo}

\pf  We proceed to prove this theorem and distinguish the following two cases.

\textbf{Case 1.} $\mathcal{N}_T\neq\emptyset$.

By Theorem \ref{3.8}, $|\mathcal{A}_T|\geq3$. It follows that there are at most $3m-4$ flexible vertices in $T$. By Lemma \ref{4.2}, in this case $T$ has at most $2\cdot 3^{m-2}$ maximum dissociation sets.

\textbf{Case 2.} $\mathcal{N}_T=\emptyset$.

We distinguish the following three subcases.

\textbf{Subcase 2.1.} $\mathcal{A}_T=\emptyset$.

By Theorem \ref{3.6}, there is no insulated $\alpha_3$-critical edge. Thus, $T$ contains exactly $m$ critical 3-paths.

If a tree $T$ is in the family $S^*_{P_3,T_1,\ldots,T_{m-1}}$, where $T_i\cong P_4 (1\leq i\leq m-1)$ (see Figure \ref{12}), then $T$ has exactly $m$ critical 3-paths and $3^{m-1}+m+1$ maximum dissociation sets. Thus the upper bound is achieved by this family of trees pictured in Figure \ref{12}. Now we show that the upper bound is achieved only on this family in this subcase.

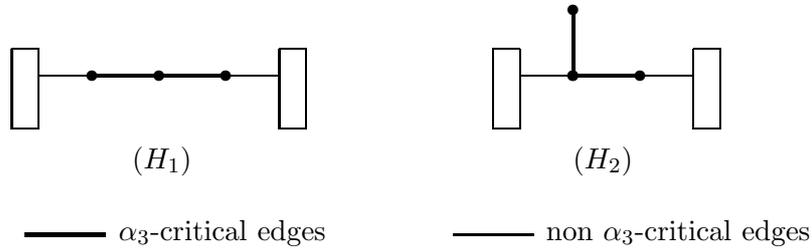
\begin{figure}[h]
	\begin{center}
		\begin{picture}(400,80)
\linethickness{0.4mm}
		\put(85,60){\line(1,0){50}}
		\thinlines
		\put(55,70){\line(1,0){10}}
		\put(55,70){\line(0,-1){30}}
		\put(65,70){\line(0,-1){30}}
		\put(55,40){\line(1,0){10}}
		\put(65,60){\line(1,0){20}}
        \put(135,60){\line(1,0){20}}
		\put(155,70){\line(1,0){10}}
		\put(155,70){\line(0,-1){30}}
		\put(165,70){\line(0,-1){30}}
		\put(155,40){\line(1,0){10}}
		\put(265,25){($H_2$)}
		\multiput(85,60)(25,0){3}{\circle*{3.5}}
		\put(100,25){($H_1$)}
		\linethickness{0.4mm}
		\put(60,0){\line(1,0){30}}
        \put(95,-2){$\alpha_3$-critical edges}
		\thinlines
		\put(220,0){\line(1,0){30}}
		\put(255,-2){non $\alpha_3$-critical edges}
		\linethickness{0.4mm}
		\put(265,60){\line(1,0){25}}
        \put(265,60){\line(0,1){25}}
		\thinlines
		\put(235,70){\line(1,0){10}}
		\put(235,70){\line(0,-1){30}}
		\put(245,70){\line(0,-1){30}}
		\put(235,40){\line(1,0){10}}
		\put(245,60){\line(1,0){20}}
        \put(290,60){\line(1,0){20}}
		\put(310,70){\line(1,0){10}}
		\put(310,70){\line(0,-1){30}}
		\put(320,70){\line(0,-1){30}}
		\put(310,40){\line(1,0){10}}
		\multiput(265,60)(25,0){2}{\circle*{3.5}}
        \multiput(265,60)(0,25){2}{\circle*{3.5}}
		\linethickness{0.4mm}
		\caption{\label{13}\small The families $H_1$ and $H_2$.}
		\end{picture}
	\end{center}
\end{figure}

\noindent\textbf{Claim 1.} Let $T$ be a tree with $3m$ vertices and $m$ critical 3-paths. If $T$ has the most maximum dissociation sets, then $T$ does not contain a subgraph that is isomorphic to a tree in the family $H_1$, or to a tree in the family $H_2$. The two families of trees $H_1$ and $H_2$ are pictured in Figure \ref{13}.

\noindent\textbf{Proof of Claim 1.} A tree in the family $H_1$ or $H_2$ has at most 10 maximum dissociation sets. When $m=3$, $10<3^{m-1}+m+1$ and any tree in the family $H_1$ or $H_2$ is not the tree with the most maximum dissociation sets.

Consider the case when $m>3$. Suppose, for a contradiction, that $T$ contains a subgraph $R$ that is isomorphic to a tree in the family $H_1$, or to a tree in the family $H_2$. Because $m>3$, the subgraph $R$ is adjacent to another critical 3-path of $T$. The larger subgraph that contains four critical 3-paths is denoted by $Q$. By simple calculation, $Q$ has at most 24 maximum dissociation sets. Thus, $T$ has at most $24\cdot 3^{m-4}$ maximum dissociation sets. This leads to a contradiction since $24\cdot 3^{m-4}<3^{m-1}+m+1$. The proof is complete.

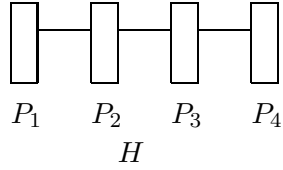
\begin{figure}[h]
	\begin{center}
		\begin{picture}(400,50)
		\put(150,45){\line(1,0){10}}
		\put(150,45){\line(0,-1){30}}
		\put(160,45){\line(0,-1){30}}
		\put(150,15){\line(1,0){10}}
        \put(180,45){\line(1,0){10}}
		\put(180,45){\line(0,-1){30}}
		\put(190,45){\line(0,-1){30}}
		\put(180,15){\line(1,0){10}}
        \put(210,45){\line(1,0){10}}
		\put(210,45){\line(0,-1){30}}
		\put(220,45){\line(0,-1){30}}
		\put(210,15){\line(1,0){10}}
        \put(240,45){\line(1,0){10}}
		\put(240,45){\line(0,-1){30}}
		\put(250,45){\line(0,-1){30}}
		\put(240,15){\line(1,0){10}}
        \put(160,35){\line(1,0){20}}
        \put(190,35){\line(1,0){20}}
        \put(220,35){\line(1,0){20}}
		\caption{\label{14}\small The family $H$.}
        \put(150,0){$P_1$}\put(180,0){$P_2$}\put(210,0){$P_3$}\put(240,0){$P_4$}
        \put(190,-15){$H$}
		\end{picture}
	\end{center}
\end{figure}

\noindent \textbf{Claim 2.} Let $T$ be a tree with $3m$ vertices and $m$ critical 3-paths. If $T$ is a tree with the most maximum dissociation sets, then there are not four critical 3-paths $P_1$, $P_2$, $P_3$ and $P_4$ such that $P_{i}$ is adjacent to $P_{i+1}$ for each $1\leq i\leq 3$, in other words, $T$ does not contain a subgraph that is isomorphic to a tree in the family $H$ pictured in Figure \ref{14}.

\noindent\textbf{Proof of Claim 2.} A tree in the family $H$ pictured in Figure \ref{14} has at most 28 maximum dissociation sets. When $m=4$, since $28<3^{m-1}+m+1$, any tree in the family $H$ is not the tree with the most maximum dissociation sets.

Consider the case when $m>4$. Suppose, for a contradiction, that $T$ contains a subgraph $R$ that is isomorphic to a tree in the families $H$. As $m>4$, the subgraph $R$ is adjacent to another critical 3-path of $T$. The larger subgraph that contains five critical 3-paths is denoted by $Q$. By calculation, $Q$ has at most 68 maximum dissociation sets. Thus, $T$ has at most $68\cdot 3^{m-5}$ maximum dissociation sets. This leads to a contradiction since $68\cdot 3^{m-5}<3^{m-1}+m+1$. The proof is complete.

By Claim 2, there exists a critical 3-path $P$ in $T$ such that every other critical 3-path is adjacent to the path $P$. By Claim 1, all other critical 3-paths are adjacent to a common vertex of the path $P$. Furthermore, if $T$ is not isomorphic to a tree in the family pictured in Figure \ref{12}, then $T$ has at most $3^{m-1}+m$ maximum dissociation sets by simple calculations.

Now we have proved that in this subcase the upper bound is achieved only on the family of trees pictured in Figure \ref{12}.

\textbf{Subcase 2.2.} $|\mathcal{A}_T|=2$.

By Theorem \ref{3.6}, in this subcase there is exactly one insulated $\alpha_3$-critical edge in $T$. Thus, each of the remaining vertices is a vertex of a critical 3-path. This contradicts the fact that the number of the remaining vertices is $3m-4$. It follows that this subcase is impossible to happen.

\textbf{Subcase 2.3.} $|\mathcal{A}_T|\geq3$.

There are at most $3m-3$ flexible vertices in $T$ in this subcase. By Lemma \ref{4.2}, $T$ has at most $3^{m-1}$ maximum dissociation sets.

Consequently, we infer that a tree $T$ on $3m$ vertices has at most ${3^{m-1}}+m+1$ maximum dissociation sets. The upper bound is achieved only in Subcase 2.1 and only on the family of trees pictured in Figure \ref{12}.
\qed

\noindent{\bf Acknowledgment.}\ We would like to thank two anonymous referees and the editor for their careful reading and helpful comments, and thank one of the referees for bringing Reference \cite{Griggs1988} to our attention.
This work was supported by the National Natural Science Foundation of China (No. 11922112), Natural Science Foundation of Tianjin,  and the Fundamental Research Funds for the Central Universities, Nankai University.

\bibliographystyle{unsrt}

\end{document}